\documentclass[10pt]{amsart}
\newtheorem{thm}{Theorem}[section]
\newtheorem{cor}[thm]{Corollary}
\newtheorem{lem}[thm]{Lemma}
\newtheorem{prop}[thm]{Proposition}

\theoremstyle{definition}
\newtheorem{defn}[thm]{Definition}
\theoremstyle{remark}
\newtheorem{rem}[thm]{Remark}
\numberwithin{equation}{section}

\def\L2{L^{2}}

\def\E{\mathcal{E}}
\def\B{\mathcal{B}}

\def\m1{^{-1}}

\def\Fr{\bigl(K , \S , \{F_i :i\in \S\}\bigr)}
\def\H{\mathcal{H}}
\def\F{\mathcal{F}}
\def\S{\mathcal{S}}

\begin{document}

\title[]{Fredholm Modules on P.C.F. Self-Similar Fractals and their Conformal Geometry.}%
\author{Fabio Cipriani, Jean-Luc Sauvageot}%
\address{Dipartimento di Matematica, Politecnico di Milano,
piazza Leonardo da Vinci 32, 20133 Milano, Italy.}%
\email{fabcip@mate.polimi.it}%
\address{Institut de Math\'ematiques, CNRS-Universit\'e Pierre et Marie
Curie, boite 191, 4 place Jussieu, F-75252 Paris Cedex 05}
\email{jlsauva@math.jussieu.fr}
\thanks{}

\date{04 July 2007}%

\begin{abstract}
The aim of the present work is to show how, using the differential
calculus associated to Dirichlet forms, it is possible to
construct Fredholm modules on post critically finite fractals by
regular harmonic structures $(D,{\bf r})$. The modules are $(d_S
,\infty )$--summable, the summability exponent $d_S$ coinciding
with the spectral dimension of the generalized laplacian operator
associated with $(D,{\bf r})$. The characteristic tools of the
noncommutative infinitesimal calculus allow to define a
$d_S$-energy functional which is shown to be a self-similar
conformal invariant.
\end{abstract}

\maketitle
\section{Introduction and description of the  results.}

The construction of Fredholm modules $(F, \H)$ on compact
topological spaces $K$ is a generalization of the theory of
elliptic differential operators on compact manifolds. In its {\it
odd} form, one requires that the elements $f$ of the algebra of
continuous functions $C(K)$ are represented as bounded operators
$\pi (f)$ on a Hilbert space $\H$ on which it is considered a
distinguished self-adjoint operator $F$ of square one $F^2 =1$,
the {\it symmetry}, in such a way that the commutators $[F,\pi
(f)]$ are compact operators.

\vskip0.2truecm\noindent This notion, introduced by M. Atiyah
[At], A.S. Mishchenko [Mis], Brown-Douglas-Fillmore [BDF], and G.
Kasparov [Kas], lyes at the core of the theory on noncommutative
differential geometry created by A. Connes [C2], where the
operator $df:=i[F,\pi (f)]$ is the operator theoretical substitute
for the differential of $f$. In its simplest example, $F$ is the
Hilbert transform acting on the space of square integrable
functions on the circle. In Atiyah's motivating case, $\H$ is the
module of square integrable sections of a smooth vector bundle
$\xi$ over a smooth manifold, the continuous functions acting in
the natural way, while $F$ arises from the parametrix of an
elliptic pseudo-differential operator of order $0$ on $\xi$.

\vskip0.2truecm\noindent In the present work, we construct Fredholm
modules on a class of {\it self-similar fractal spaces}, known as
{\it post critically finite} (shortened as p.c.f. since now on).
Self-similarity refers to the fact that such a space can be
reconstructed as finite union of homeomorphic pieces of itself.
The p.c.f. property, on the other hand, translates or generalizes
mathematically, a property of {\it finite ramification} and it is
for this reason that, in general, these spaces fail to be
manifolds modelled on open Euclidean sets, so that the usual
Leibniz-Newton infinitesimal calculus is no more avaible.

\vskip0.2truecm\noindent These spaces have been largely
investigated from the point of view of potential and spectral
analysis (Dirichlet forms, Laplacians, heat kernels, Green
functions, eigenvalues distribution) and probability theory
(construction and analysis of diffusive Markov processes) (see for
example [Ba], [FS], [Ki], [Ku]). Spaces of this class, including,
for example, Koch's curve, Sierpinski's gasket, Hata's tree-like
set and Lindstr$\o$m's snowflake, exhibit singular behaviors when
compared, from the above points of view, to differentiable
riemannian manifolds. For example: i) their energy measures are, in
general, singular with respect to any self-similar volume measure
(see [BST], [Hi]); ii) in the strong symmetric case, they support
localized eigenfunctions (see [FS], [Ki]); iii) in the so called
{\it arithmetic or lattice} case, the integrated density of states
is discontinuous (see [FS], [Ki]). It is worth to mention that the
study of these exotic behaviors of fractals spaces was suggested
and motivated for application to condensed matter physics (see
[L], [RT]).

\vskip0.2truecm\noindent The first constructions of Fredholm
modules over subsets of nonintegral Hausdorff dimension were given
by A. Connes in [C2 IV.3] for quasi-circles embedded in the plane
and Cantor subsets of the real line, while D. Guido - T. Isola
considered in [GI1,2,3] more general subsets of $\mathbb{R}^n$.
\par\noindent
Our construction of Fredholm modules on p.c.f. fractals is based
on the notion of {\it regular harmonic structure} introduced by J.
Kigami [Ki] and on the differential calculus associated to
Dirichlet forms we developed in [CS].
\par\noindent
To any fixed harmonic structure on a p.c.f fractal one can
associate its self-similar Dirichlet form $(\E ,\F )$. This is a
lower semi-continuous quadratic form, defined on a uniformly dense
subalgebra $\F$ of $C(K)$ and satisfying the characteristic
contraction property
\[
\E [a\wedge 1]\le \E[a]
\]
which generalizes
the Dirichlet integral of an Euclidean space.\par\noindent Dirichlet forms
can be canonically represented as graphs semi-norms
$\E[a]=\|\partial a\|_{\H}^2$ of an essentially unique {\it
derivation} $\partial : \F\rightarrow \H$ [CS]. This is a map, taking
values in a Hilbert space $\H$ which is a module over the algebra
$C(K)$, and satisfying the Leibnitz rule:
\[
\partial (ab)=(\partial
a)\cdot b + a\cdot (\partial b)\, .
\]
It is by the derivation
$\partial$ that the Dirichlet form $\E$ defines, in a natural way,
 a differential calculus on the fractal $K$.\par\noindent
 To quantize this calculus, we then
define the Fredholm module $(F,\H )$ by the symmetry $F$
corresponding to the subspace ${\it Im}\,\partial\subset \H$: $F$
acts as the identity on the range ${\it Im}\partial$ of the
derivation and specularly on its orthogonal complement $({\it
Im}\,\partial )^\perp$.\vskip0.2truecm\noindent It is worth to
recall that the Dirichlet form $\E$ is a quadratic form closable
on the Lebesgue space $L^2 (K,\mu )$ with respect to a large set
of positive Radon measures $\mu$ (see [Ki]). By classical results
of Dirichlet form theory (see [BD], [FOT]), the closure of $\E$
with respect to such a measure $\mu$ is then the quadratic form of
a positive, self-adjoint operator $\Delta_\mu$ on $L^2 (K,\mu )$
which generates a Markovian semigroup $e^{-t\Delta_\mu}$ whose
heat kernel $p(t,x,y)$ gives the transition probabilities of a
Markovian diffusion process on $K$.


\vskip0.2truecm\noindent The above definition
for $(\F ,\H)$ is inspired by a result of Connes-Sullivan [Co
IV.4] concerning a canonical construction of a Fredholm module on
a even dimensional manifold $V$. Their construction makes use of a
fixed riemannian metric on $V$ but the resulting Fredholm module
directly determines the underlying conformal structure of $V$.
This is explicitly seen by a formula reproducing the fundamental
conformal invariant, namely the ${\rm dim}V$-homogeneous Dirichlet
integral $\int_V |\nabla a|^{{\rm dim}V}$ through a suitable
summation procedure known as the {\it Dixmier trace}.\par\noindent
In our setting, by analyzing the speed of vanishing of the
sequences of the eigenvalues of the commutators $[F,\pi (f)]$,
through the use of the Schatten's classes of compact operators and
other interpolation ideals, we are able to construct, still
through the Dixmier trace, a new, densely defined, strongly local,
convex energy functional $\E_C$ on $C(K)$. Its homogeneity
exponent $d_S$ equals the {\it spectral dimension} of the heat
semigroup generated by the generalized Laplacian associated to the
closure of the Dirichlet form with respect to a natural
self-similar measure on K, called by J. Kigami and M. Lapidus [KL]
the {\it Riemannian volume measure} of $K$ (at least for decimable
fractals). Showing that $\E_C$ is self-similarly invariant one
could be inclined to consider $\E_C$ as defining a {\it
generalized conformal structure} on $K$. \vskip0.2truecm\noindent
We finally remark that our construction allows to associate to
each harmonic structure a topological invariant of $K$, namely the
K-homology class of the Fredholm module $(\F ,\H)$ (see [BDF]).
Using the Chern character in cyclic cohomology [C2] it should be
challenging to measure how harmonic structures on the same p.c.f.
fractal may differ from a topological point, of view so to
consider a moduli space of harmonic structures.

\section{Laplacian and Dirichlet forms on P.C.F. Self-Similar Sets}
In this section we will briefly recall, for reader's convenience,
the main definitions and properties of the objects we will
investigate. See [Ki] for details.

\begin{defn}{\bf (Self-similar structures)}
Let $K$ be a compact metrizable topological space and let $\S:=
\{1,2,\dots ,N\}$ for a fixed integer $N$ greater than one. For
each $i\in \S$, let us denote by $F_i$ a fixed continuous injection
of $K$ into itself. Then, $\bigl(K , \S , \{F_i :i\in \S\}\bigr)$
is called a {\it self-similar structure} if there exists a
continuous surjection $\pi : \Sigma \rightarrow K$ such that $F_i
\circ \pi = \pi \circ \sigma_i$ for every $i\in \S$, where $\Sigma
:=\S^\mathbb{N}$ is the one-sided shift space and
$\sigma_i:\sum\rightarrow\sum$ denotes the injection $\sigma_i (w_1 w_2 w_3 \dots
):=iw_1 w_2 w_3 \dots $ for each $w_1 w_2 w_3 \dots \in\sum$.
\end{defn}

Notice that if $\Fr$ is a self-similar structure, then $K$ is
self-similar in the sense that

\begin{equation}
K=\bigcup_{i\in\S} F_i (K)\, .
\end{equation}

It is customary to denote by $W_m =: \S^m$ the set of words of
length $m\in \mathbb{N}$, composed using the letters of the
alphabet $\S$, with the understanding that $W_0 := \emptyset$,
setting also $W_* := \bigcup_{m\in \mathbb{N}} W_m$ for the whole
vocabulary. Each word $w=w_1\dots\ w_m\in W_m$ defines a
continuous injection $F_w :K\rightarrow K$ by $F_w :=
F_{w_1}\circ\dots\circ F_{w_m}$, whose image $F_w (K)$ is denoted
by $K_w$.

\begin{defn}{\bf (Post critically finite fractals)}
Let $\bigl(K , \S , \{F_i :i\in \S\}\bigr)$ be a self-similar
structure. The {\it critical} set $\mathcal{C}\subset \Sigma$ and
the {\it post critical} set $\mathcal{P}\subset \Sigma$ are
defined by
\[
\mathcal{C}:= \pi^{-1} \Bigl( \bigcup_{i\ne j} K_i\cap
K_j\Bigr)\qquad {\rm and}\qquad \mathcal{P}:= \bigcup_{n\ge 1}
\sigma^n (\mathcal{C})\, ,
\]
where $\sigma :\Sigma\rightarrow\Sigma$ is the shift map defined
by $\sigma (w_1 w_2\dots ):= w_2 w_3\dots$. A self-similar
structure is called {\it post critically finite} (p.c.f. for
short) provided $\mathcal{P}$ is a finite set. One sets also $V_0
:= \pi (\mathcal{P})$, considered as the boundary of K, and
\[
V_m := \bigcup_{w\in W_m} F_w (V_0)\qquad{\rm and}\qquad V_*
:=\bigcup_{m\in \mathbb{N}}V_m\, .
\]
Its is easy to see that $V_m \subset V_{m+1}$ and that $V_*$ is
dense in $K$.
\end{defn}
For a finite set $V$, equipped with the standard counting measure,
denote by $l(V)$ the space of scalar functions on $V$ with its
scalar product $(u|v):=\sum_{p\in V} u(p)v(p)$ and by
$\mathcal{L}(V)$ the collection of the {\it Laplacian operators}
on $V$, i.e. the generators $
L$ of the conservative, {\it symmetric Markovian
semigroups} $e^{-tL}$ on $l(V)$. These are, essentially, symmetric, positive definite matrices $L=\{L_{u,v} : u,v\in V\}$ such that  $L_{u,v}\le 0$ for $u\ne v$ and $\sum_{u\in V} L_{u,v} =0$ for all $v\in V$.\par\noindent
 For $L\in \mathcal{L}(V)$ let
$\E_L$ the associated Dirichlet form on $l(V)$: $\E_L
(u,v)=(Lu|v)$.

For a fixed self-similar structure $\bigl(K , \S , \{F_i :i\in
\S\}\bigr)$ on $K$, a Laplacian $D\in \mathcal{L}(V_0)$ and a
vector $\mathbf{r} := (r_1 ,\dots ,r_N)$, where $r_i >0$ for $i\in
\S$, define for each $m\ge 0$ the quadratic form
\begin{equation}
\E^{(m)} [u]:=\sum_{w\in W_m} \frac{1}{r_w} \E_D [u\circ F_w]\,
,\qquad\qquad u\in l(V_m)
\end{equation}
where $r_w := r_{w_1}\dots r_{w_m}$ for $w=w_1\dots w_m\in W_m$.
It is easy to see that there exists $H_m\in \mathcal{L}(V)$ such
that $\E^{(m)} [u]=(H_m u |u)$.
\vskip0.2truecm\noindent
We now introduce the main object of analysis on fractals.
\begin{defn}{\bf (Harmonic structures)}
$(D,\mathbf{r})$ is said to be a {\it harmonic structure} on
$\bigl(K , \S , \{F_i :i\in \S\}\bigr)$ if for all $m\ge 0$ and
for any $u\in l(V_m)$ one has
\begin{equation}
\E^{(m)} [u]= \min\{\E^{(m+1)} [v]:v\in l(V_{m+1})\, ,\,\,\,
v|_{V_m} =u\}.
\end{equation}
\end{defn}
It is known that (2.3) holds for all $m\ge 0$ if and only if it
holds for $m=0$.
\begin{defn}
If $(D,\mathbf{r})$ is a {\it harmonic structure} on $\bigl(K , \S
, \{F_i :i\in \S\}\bigr)$, define
\begin{equation}
\F:= \{u\in l(V_*): \lim_{m\to \infty}
\E^{(m)}[u|_{V_m}]<\infty\}\, ,\qquad \F_0:= \{u\in \F: u|_{V_0}
=0\}
\end{equation}
and
\begin{equation}
\E [u]:= \lim_{m\to \infty} \E^{(m)}[u|_{V_m}]\qquad {\rm
for}\qquad u\in \F\, .
\end{equation}
\end{defn}
Since the quadratic form $(\E, \F )$ is defined in a self-similar
way, it naturally satisfies the following self-similarity.
\begin{prop}{\bf [Ki] (Self-similar quadratic form)} Let $(D,\mathbf{r})$ be a {\it harmonic structure} on
$\bigl(K , \S , \{F_i :i\in \S\}\bigr)$. Then
$u\in \F$ if and only if $u\circ F_i\in \F$ for all $i\in \S$ and
in that case
\begin{equation}
\E [u]=\sum_{i\in \S} \frac{1}{r_i} \E [u\circ F_i]\,
,\qquad\qquad u\in \F .
\end{equation}
\end{prop}
To construct Dirichlet forms on $K$ we need to fix measures on
$K$. Here is a natural class one may consider.
\begin{prop}{\bf [Ki] (Self-similar measure)}
For a fixed vector of weights $(\mu_1 ,\dots ,\mu_N)$ with $\mu_i
>0$ for $i\in \S$ and $\sum_{i\in \S} \mu_i =1$, there exists a
unique Borel measure $\mu$ on $K$ such that
\begin{equation}
\int_K f d\mu = \sum_{i\in\S}\mu_i \int_K f\circ F_i\, d\mu\,
,\qquad\qquad f\in C(K).
\end{equation}
$\mu$ is called the
self-similar measure with weights $(\mu_1 ,\dots ,\mu_N).$
\end{prop}
\begin{thm}{\bf [Ki] (Dirichlet forms and generalized laplacians)}
Let $(D,\mathbf{r})$ be a {\it harmonic structure} on a p.c.f.
self-similar structure $\bigl(K , \S , \{F_i :i\in \S\}\bigr)$,
$\mu$ the self-similar measure on $K$ with weights $(\mu_1 ,\dots
,\mu_N)$ and assume that $\mu_i r_i <1$ for all $i\in
\S$.\par\noindent Then $\F$ is naturally embedded in $L^2
(K,\mu)$, $(\E ,\F)$ and $(\E ,\F_0)$ are regular, local Dirichlet
form on $K$ and their associated nonnegative self-adjoint
operators $H_N$ and $H_D$ have compact resolvent.
\end{thm}


\begin{defn}{\bf (Eigenvalues distribution)}
Assuming the same hypotheses as in the previous theorem,
define, for $*=N,D$, the eigenspace corresponding to $\lambda\in
\mathbb{R}$ as
\begin{equation}
E_* (\lambda )=\{u\in D(H_* ): H_* u = \lambda u\}\, .
\end{equation}
If the multiplicity ${\rm dim}\, E_* (\lambda )$ is not zero then
$\lambda$ is said to be an $*$-eigenvalue and a non zero $u\in E_*
(\lambda )$ is said to be a $*$-eigenfunction belonging to the
$*$-eigenvalue $\lambda$. The collection ${\rm Sp} (H_*)$ of all
the eigenvalues of $H_*$ is called the spectrum of $H_*$. As $H_*$
is unbounded with compact resolvent, ${\rm Sp} (H_*)$ is unbounded
and discrete, consisting of isolated eigenvalues of finite
multiplicity only.\par\noindent The function
\begin{equation}
\rho_* (\cdot ,\mu ): \mathbb{R}\rightarrow \mathbb{N}\qquad
\rho_* (x ,\mu ) := \sum_{ \lambda\le x} {\rm dim}\, E_* (\lambda
)
\end{equation}
is called the {\it eigenvalues counting function} of $H_*$. As
$H_*$ is nonnegative and unbounded $\rho_* (x ,\mu )=0$ if $x<0$
and $\lim_{x\to +\infty}\rho_* (x ,\mu )=+\infty$.
\end{defn}
The following is the fractal analogue of the famous Weyl's
asymptotic formula for the eigenvalue distribution of the laplacian on a compact riemannian manifold.
\begin{thm}{\bf [Ki]}
Let $(D,\mathbf{r})$ be a {\it harmonic structure} on a p.c.f.
self-similar structure $\bigl(K , \S , \{F_i :i\in \S\}\bigr)$,
$\mu$ the self-similar measure on $K$ with weights $(\mu_1 ,\dots
,\mu_N)$ and assume that $\mu_i r_i <1$ for all $i\in \S$.
\par\noindent Let $d_S = d_S(\mu )$ be the unique positive real
number satisfying
\begin{equation}
\sum_{i\in \S} \gamma_i^{d_S} = 1\, ,
\end{equation}
where $\gamma_i := \sqrt{r_i\mu_i}$ for $i\in \S$. $d_S$ is called the spectral exponent of $(\E\, ,\F\, ,\mu )$. Then
\begin{equation}
0<\liminf_{x\to +\infty} \frac{\rho_* (x ,\mu )}{x^{d_S/2}}\le
\limsup_{x\to +\infty} \frac{\rho_* (x ,\mu )}{x^{d_S/2}}
<+\infty\, ,
\end{equation}
the $\liminf$ and the $\limsup$ are the same for $*=N$ and $*=D$. In the {\it
non-lattice case}, where $\sum_{i\in \S} \mathbb{Z}\log \gamma_i$
is dense subgroup of $\mathbb{R}$, defining
\[
R(x):=\rho_D (x,\mu ) - \sum_{i\in\S} \rho_D (\gamma_i^2x,\mu
)\qquad U(t):=e^{-td_\S} R (e^{2t})
\]
 we have
\begin{equation}
\lim_{x\to +\infty} \frac{\rho_* (x ,\mu )}{x^{d_S/2}} =
\Bigl(-\sum_{i\in\S}\gamma_i^{d_S}\log\gamma_i^{d_S}\Bigr)^{-1}
d_\S\int_\mathbb{R} U(t)dt\, .
\end{equation}
\end{thm}
Comparing the above result with the Weyl's classical one for
compact remannian manifold, one is led to define
\begin{defn}
Let $(D,\mathbf{r})$ be a {\it harmonic structure} on a p.c.f.
self-similar structure $\bigl(K , \S , \{F_i :i\in \S\}\bigr)$,
$\mu$ the self-similar measure on $K$ with weights $(\mu_1 ,\dots
,\mu_N)$ and assume that $\mu_i r_i <1$ for all $i\in \S$. The
{\it spectral volume} ${\rm vol} (K,\mu)$ is then defined by
\begin{equation}
{\rm vol} (K,\mu) :=
\Bigl(-\sum_{i\in\S}\gamma_i^{d_S}\log\gamma_i^{d_S}\Bigr)^{-1}
d_\S\int_\mathbb{R} U(t)dt\, .
\end{equation}
\end{defn}
On compact riemannian manifolds, the Connes' trace formula [C1] allows to reconstruct
the riemannian measure through the knowledge of suitable eigenvalues distributions. This is done by the {\it Dixmier trace}
${\rm Tr}_\omega$, a trace functional on the space of compact operators on a Hilbert space, depending on the choice
of certain ultrafilters on $\mathbb{R}_+$. This functional is singular
in the sense that it vanishes on the ideal of trace-class operators. The following generalization of the Connes' trace formula
has been proved on fractals by J. Kigami and M. Lapidus.
\begin{thm}{\bf [KL]}
Let $(D,\mathbf{r})$ be a {\it harmonic structure} on a p.c.f.
self-similar structure $\bigl(K , \S , \{F_i :i\in \S\}\bigr)$,
$\mu$ the self-similar measure on $K$ with weights $(\mu_1 ,\dots
,\mu_N)$ and assume that $\mu_i r_i <1$ for all $i\in \S$. Then
there exists a unique positive Borel measure $\nu_\mu$ on $K$ such
that
\begin{equation}
\int_K fd\nu_\mu = {\rm Tr}_\omega \Bigl(f\circ H_D^{-d_\S
/2}\Bigr)\, ,
\end{equation}
where the symbol $f$ denotes both a continuous function on $K$ as
the associated multiplication operator on $L^2 (K,\mu)$. Moreover
the total mass of $\nu_\mu$ equals the spectral volume of $K$:
$\nu_\mu (K)={\rm vol} (K,\mu)$.
\end{thm}
It has been proved in [KL], that for certain classes of fractals,
$\nu_\mu$ is the self-similar measure on $K$ with weights $\nu_i =
\gamma_i^{d_\S}$.


\section{Fredholm Modules associated to Harmonic Structures on p.c.f. fractals.}
In this section we consider a fixed harmonic
structure $(D,\mathbf{r})$ on a p.c.f. self-similar structure
$\bigl(K , \S , \{F_i :i\in \S\}\bigr)$.
\par\noindent
Choosing a self-similar measure $\mu$ on $K$ with weights $(\mu_1
,\dots ,\mu_N)$ such that $\mu_i r_i <1$ for all $i\in \S$, we can
consider, by Theorem 2.7, the Dirichlet form $(\E ,\F )$
associated to $(D,\mathbf{r})$ on $L^2 (K,\mu)$.
\par\noindent Applying the general theory developed in [CS], it is possible to consider
a {\it differential calculus on the fractal} $K$, associated to the Dirichlet form $(\E ,\F )$.
In other words:
\begin{prop}{\bf [CS]}
There exists an essentially unique {\it derivation} $\partial : \B\rightarrow
\H$, defined on the {\it Dirichlet algebra} $\B = C(K)\cap \F$
with values in a {\it real Hilbert module} $\H$, which is a {\it differential square root}
of the Dirichlet form in the precise sense that
\begin{equation}
\E [u]=\|\partial u\|^2_\H\qquad u\in \B\, .
\end{equation}
By this, we mean that $\H$ is a Hilbert space on which the algebra $C(K)$ acts continuously in such
a way that the {\it Leibniz rule} holds true:
\begin{equation}
\partial (ab) = (\partial a)b + a(\partial b)\, \qquad a,b\in \B\, .
\end{equation}
In turn, the self-adjoint operator $H_N$ associated to $(\E ,\F)$ on $L^2 (K,\mu)$ appears as a {\it generalized laplacian}
\begin{equation}
H_N = \partial_\mu ^*\circ \partial_\mu
\end{equation}
where $\partial_\mu$ denotes the closure of $(\partial ,\B)$ in $L^2 (K,\mu)$. A corresponding result holds true for $H_D$.\par\noindent
\end{prop}

When the harmonic structure is regular, i.e. $r_i<1$ for $i\in\S$,
then the Dirichlet algebra $\B$ coincide with the domain $\F$ of
the Dirichlet form, which, in this case, is already a sub-algebra
of $C(K)$ [Ki 3.3].
\vskip0.2truecm
To recover the information potentially carried by the derivation, we consider its associated {\it phase}.
\begin{defn}
Let us consider a fixed harmonic structure $(D,\mathbf{r})$ on a
p.c.f. self-similar structure $\bigl(K , \S , \{F_i :i\in
\S\}\bigr)$. Let $\partial :\B\rightarrow \H$ be the associated
derivation, defined on the Dirichlet algebra $\B = C(K)\cap \F$
with values in the symmetric Hilbert module $\H$. Let $P\in {\rm
Proj}\, (\H)$ the projection onto the closure $\overline{{\rm
Im}\partial}$ of the range of the derivation
\begin{equation}
P\H = \overline{{\rm Im}\partial}
\end{equation}
and $F=P-P^\perp :\H\rightarrow\H$ the associated {\it phase} or
{\it symmetry}.
\end{defn}
The following result shows that the phase operator associated to a regular harmonic structure on p.c.f.
fractals is an elliptic operator on $K$ in the sense of M. Atiyah [At].
\begin{thm}{\bf (Fredholm module structure on fractals)}
Let $(D,\mathbf{r})$ be a fixed regular harmonic structure on a
p.c.f. self-similar structure $\bigl(K , \S , \{F_i :i\in
\S\}\bigr)$. Then $(F,\H)$ is a {Fredholm module over $C(K)$} in
the sense of [At] and a {\it densely 2-summable Fredholm module
over $C(K)$} in the sense of [C IV 1.$\gamma$ Definition 8].
\end{thm}
\begin{proof}
Clearly $F^* = F$, $F^2 = I$. Let us start to prove that $[F,a]$
is Hilbert-Schmidt for all real valued $a\in\F$. Since
\begin{equation}
[P,a] = PaP^\perp -P^\perp aP
\end{equation}
and $a$ is real valued, we have
\begin{equation}
|[P,a]|^2 = |PaP^\perp |^2 + |P^\perp aP|^2
\end{equation}
so that
\begin{equation}
\|[F,a]\|^2_{\mathcal{L}^2} = 4\|[P,a]\|^2_{\mathcal{L}^2} = 8
\|P^\perp aP\|^2_{\mathcal{L}^2}\, .
\end{equation}
Using the Leibnitz rule for the derivation $\partial $, the fact
that $P\circ
\partial =
\partial$ and $P^\perp\circ
\partial = 0$, we have, for all $b\in\F$,
\begin{equation}
P^\perp aP (\partial b)=P^\perp (a\partial b) = P^\perp (\partial
(ab)-(\partial a)b) = -P^\perp ((\partial a)b)
\end{equation}
so that
\begin{equation}
\|P^\perp aP (\partial b)\| = \|P^\perp ((\partial a)b)\| \le
\|(\partial a)b\|\, .
\end{equation}
Let us choose a self-similar Borel measure $\mu$ on $K$ with
weights $(\mu_1 \, ,\dots\, ,\mu_N)$ such that $r_i\mu_i<1$ for
$\in\S$. By [Ki Theorem 3.4.7], $(\E, \F_0 )$ has discrete
spectrum $\{0<\lambda_1 \le \lambda_2 < \dots\}$ in $L^2 (K,\mu
)$. Denoting by $a_1 ,a_2 , \dots$ the corresponding
eigenfunctions, we have that the vectors $\xi_k :=
\lambda^{-1/2}_k\partial a_k$, $k\ge 1$, form an orthonormal
complete system in $P\H$. Then by (3.7) and (3.9)
\begin{equation}
\begin{split}
\|[F,a]\|^2_{\mathcal{L}^2} &= 8 \|P^\perp aP\|^2_{\mathcal{L}^2}
= 8\sum_{k=1}^\infty \lambda^{-1}_k \|P^\perp aP(\partial
a_k)\|^2_\H \le 8 \sum_{k=1}^\infty \lambda^{-1}_k \|(\partial
a )a_k\|^2_\H \\&= 8\sum_{k=1}^\infty \lambda^{-1}_k \int_K a_k^2\,
d\Gamma (a) = 8 \int_K
\bigl(\sum_{k=1}^\infty\lambda^{-1}_k a_k^2 \bigr)d\Gamma (a)\\
&=8\int_K
g\,d\Gamma (a)\\
\end{split}
\end{equation}
where $g$ is the restriction on the diagonal of $K\times K$ of the Green function $G(x,y)=\sum_{k=1}^\infty\lambda^{-1}_k a_k (x)a_k (y)$, kernel of the compact operator $H^{-1}_{\rm D}$ on $L^2
(K,\mu )$ (see [Ki 3.6]) and $\Gamma (a)$ is the {\it energy measure} of $a\in\F$ defined by the Dirichlet form [FOT]
\[
\int_K b\, d\Gamma (a) := (\partial a |(\partial a)b) = \E (a |ab)-\frac{1}{2}\E (b\,|a^2)\, ,\qquad b\in\B=C (K) \cap \F\, .
\]
Since the harmonic structure is regular, $G$ is
continuous on $K\times K$ ([Ki Proposition 3.5.5]) and we have
\[
\|[F,a]\|^2_{\mathcal{L}^2} \le 8 \Bigl(\sup_{x\in K} g(x)\Bigr)
\E[a]<+\infty
\]
for all $a\in \F\subset C(K).$ Since $\F$ is uniformly dense in
$C(K)$, $[F,a]$ is norm continuous with respect to $a\in C(K)$ and
the space of compact operators is norm closed, we have that
$[F,a]$ is compact for all $a\in C(K)$.
\end{proof}
\rem The proof given above shows that the regularity of a function $a\in\B = \F\cap C(K)$, can be detected using
the energy form $\E$ and an auxiliary reference measure $\mu$ with respect to which $\E$ has discrete spectrum
(in this respect see [Ki Theorem 3.4.6, Corollary 3.4.7]).
The effectiveness of the upper bound on the Hilbert-Schmidt norm of the commutator $[F,a]$ depends on the integrability
of the diagonal part of the Green function (of $\E$ with respect to $\mu$) with respect to the energy measure $\Gamma (a)$.
The same proof thus provides a method for constructing Fredholm modules even for non regular harmonic structures.
In these situations, one no more has $\F\subseteq C(K)$ but may uses the core of harmonic functions associated to the harmonic structure
\[
\H_* =: \{a\in \F : {\rm a\,\, is\,\, an\,\, m-harmonic\,\, function\,\,
for\,\,  some\,\,  m\ge 0}\}\subset \B\, .
\]
In particular see [Ki 3.2] and the proof of [Ki Theorem 3.4.6].
\vskip0.2truecm\noindent
We are now interested to investigate finer summability properties of
the commutators $[F,a]$. The following lemma contains an estimate we will need below. It is essentially [Ki Lemma 5.3.5].
\begin{lem}
Let $(D,\mathbf{r})$ be a fixed regular harmonic structure on a
p.c.f. self-similar structure $\bigl(K , \S , \{F_i :i\in
\S\}\bigr)$ and let $\mu$ be a self-similar Borel measure on $K$ with
weights $(\mu_1 \, ,\dots\, ,\mu_N)$ such that $r_i\mu_i<1$ for
$i\in\S$. Denote by $d_S$ the spectral exponent of $(\E\, ,\F\, ,\mu )$. Then the potential operators
\begin{equation}
G_p := \int_0^\infty t^{\frac{p}{2}-1}e^{-tH_D} dt\, ,\qquad d_S <
p \le 2
\end{equation}
are compact in $L^2 (K,\mu)$ and their integral kernels $g_p$ are
positive continuous functions satisfying, for some $c_1>0$,
\begin{equation}
g_p (x,y) \le c_1\Bigl(\frac{1}{\lambda_1} + \frac{2}{p-d_S}
\Bigr).
\end{equation}
\end{lem}
\begin{proof}
By the Spectral Theorem $G_p = \Gamma (\frac{p}{2})
H_D^{-\frac{p}{2}}$, so that the compactness follows from the
discreetness of the spectrum of the laplacian. Let $p_D (t,x,y)$
be the kernel of the heat semigroup $e^{-tH_D}$ so that
\begin{equation}
g_p (x,y) = \int_0^\infty t^{\frac{p}{2}-1}p_D (t,x,y) dt\, .
\end{equation}
By [Ki Lemma 5.3.5] there exists $c_1>0$ such that
\begin{equation}
p_D (t,x,y) \le\begin{cases}
    c_1t^{-\frac{d_S}{2}}\qquad t\in(0,1] \\
    c_1e^{-(t-1)\lambda_1}\qquad t\ge 1\, ,
  \end{cases}
\end{equation}
(in fact $c_1=\|e^{-H_D}\|_{L^1\to L^\infty}$) from which we get
\begin{equation}
g_p (x,y) \le c_1\Bigl\{\int_0^1 t^{\frac{p}{2}-1}\cdot
t^{-\frac{d_S}{2}} + \int_1^\infty t^{\frac{p}{2}-1}\cdot
e^{-(t-1)\lambda_1}\Bigr\} \le c_1\Bigl(\frac{1}{\lambda_1} +
\frac{2}{p-d_S} \Bigr).
\end{equation}

\end{proof}

\begin{thm}{\bf (Commutators and Shatten classes)}
Let $(D,\mathbf{r})$ be a fixed regular harmonic structure on a
p.c.f. self-similar structure $\bigl(K , \S , \{F_i :i\in
\S\}\bigr)$ and let $\mu$ be a self-similar Borel measure on $K$ with
weights $(\mu_1 \, ,\dots\, ,\mu_N)$ such that $r_i\mu_i<1$ for
$\in\S$.\par\noindent Then $(F,\H)$ is a densely p-summable
Fredholm module over $C(K)$ for all $d_S < p\le 2$. In particular
\begin{equation}
{\rm Trace} \Bigl(|[F,a]|^p\Bigr) \le c_2 (p)^{\frac{p}{2}}
\cdot\Bigl(\E[a]\Bigr)^{\frac{p}{2}}\cdot\Bigl[{\rm Trace}
\Bigl(H_D^{-\frac{p}{2}}\Bigl)\Bigr]^{1-\frac{p}{2}} \, ,
\end{equation}
where $c_2 (p) := 16c_1(\frac{1}{\lambda_1} + \frac{2}{p-d_S} )
\Gamma (\frac{p}{2})^{-1}$, for all a $\in \F$ ($c_1$ being the
constant appearing in Lemma 3.5).
\end{thm}
\begin{proof}
Let us fix $a\in \F$ real valued and denote by $\{\mu_k (T): k\ge
0\}$ the non vanishing singular values of a compact operator $T$
arranged in decreasing order and repeated according to their
multiplicity. Recall that $\mu_k (T) = \mu_k (T^*)$. Setting
$S:=[P,a]$ and $T=:P^\perp aP$, from (3.5) we get
\begin{equation}
|S| = |T^*| + |T|
\end{equation}
and then
\begin{equation}
\mu_{n+m} (S) = \mu_{n+m} (|T^*| + |T|) \le \mu_{n} (T^*) +
\mu_{m} (T) = \mu_{n} (T) + \mu_{m} (T)
\end{equation}
\begin{equation}
\mu_k (S) \le 2\mu_{[\frac{k}{2}]} (T)\le 2\mu_k (T),\qquad k\ge 0
\end{equation}
and finally
\begin{equation}
\begin{split}
{\rm Trace}\, (|S|^p) &= \sum_{k=0}^\infty \mu_k (|S|^p) =
\sum_{k=0}^\infty \mu_k (S)^p \\
&\le  2^p
\sum_{k=0}^\infty \mu_k (T)^p = 2^p\,{\rm Trace}\, (|T|^p)\, .
\end{split}
\end{equation}
Since $\{\xi_k := \lambda_k^{-1/2} \partial a_k : k\ge 1\}$ is a
complete orthonormal family in the Hilbert space $P\H$ and, by
assumption, $p\le 2$, we can use inequality [S Remark 1 page 17]
and $(3.9)$ to get
\begin{equation}
\begin{split}
{\rm Trace}\, (|T|^p) &= \sum_{k=0}^\infty \mu_k (T)^p \le
\sum_{k=0}^\infty \|T\xi_k\|^p = \sum_{k=0}^\infty
\bigl(\|T\xi_k\|^2\bigr)^{p/2} \\
& = \sum_{k=0}^\infty \bigl(\|\lambda_k^{-1/2} (\partial
a)a_k \|^2\bigr)^{p/2}\\
& = \sum_{k=0}^\infty \Bigl(\int_K \lambda^{-1}_k a_k^2\, d\Gamma
(a) \Bigr)^{p/2}\, .
\end{split}
\end{equation}
By H\"older's inequality in the spaces $l^q (\mathbb{N})$, with
conjugate exponents $2/p$ and $2/(2-p)$, we obtain
\begin{equation}
\begin{split}
{\rm Trace}\, (|T|^p) &\le \sum_{k=0}^\infty \Bigl(\int_K
\lambda^{-1}_k a_k^2\, d\Gamma (a) \Bigr)^{p/2} \\
& = \sum_{k=0}^\infty \lambda^{-\frac{p(2-p)}{4}}_k\Bigl(\int_K
\lambda^{-\frac{p}{2}}_k a_k^2\, d\Gamma (a) \Bigr)^{p/2} \le \\
& = \Bigl(\sum_{k=0}^\infty
\lambda^{-\frac{p}{2}}_k\Bigr)^{1-\frac{p}{2}}\Bigl(\int_K
\sum_{k=0}^\infty \lambda^{-\frac{p}{2}}_k a_k^2\, d\Gamma (a)
\Bigr)^{p/2} \\
& = \Bigl[ {\rm Trace}\, (H_D^{-\frac{p}{2}})
\Bigr]^{1-\frac{p}{2}}\Bigl(\int_K \sum_{k=0}^\infty
\lambda^{-\frac{p}{2}}_k a_k^2\, d\Gamma (a) \Bigr)^{p/2}\, .
\end{split}
\end{equation}
Since $G_p = \Gamma (p/2) H_D^{-\frac{p}{2}}$, we have
$\sum_{k=0}^\infty \lambda^{-\frac{p}{2}}_k a_k^2 (x) = \Gamma
(p/2)^{-1} g_p (x,x)$. From Lemma 3.5 and (3.22) we have
\begin{equation}
\begin{split}
{\rm Trace}\, (|T|^p) &\le \Gamma (p/2)^{-\frac{p}{2}}\Bigl[ {\rm
Trace}\, (H_D^{-\frac{p}{2}}) \Bigr]^{1-\frac{p}{2}}\Bigl(\int_K
g_p
(x,x)\, \Gamma (a)(dx) \Bigr)^{p/2} \\
& \le c_1^{\frac{p}{2}}\Bigl(\frac{1}{\lambda_1} + \frac{2}{p-d_S}
\Bigl)^{\frac{p}{2}} \Gamma
(p/2)^{-\frac{p}{2}}\Bigl(\E[a]\Bigr)^{p/2}\Bigl[ {\rm Trace}\,
(H_D^{-\frac{p}{2}}) \Bigr]^{1-\frac{p}{2}}\, .
\end{split}
\end{equation}
Noticing that $[F,a]=2S$, we finally obtain from (3.20) and (3.23)
\begin{equation}
\begin{split}
{\rm Trace} \Bigl(|[F,a]|^p\Bigr) & \le 4^p
c_1^{\frac{p}{2}}\Bigl(\frac{1}{\lambda_1} + \frac{2}{p-d_S}
\Bigl)^{\frac{p}{2}} \Gamma
(p/2)^{-\frac{p}{2}}\Bigl(\E[a]\Bigr)^{p/2}\Bigl[ {\rm
Trace}\, (H_D^{-\frac{p}{2}}) \Bigr]^{1-\frac{p}{2}} \\
& = c_2 (p)^{\frac{p}{2}}\cdot \Bigl(\E[a]\Bigr)^{p/2}\Bigl[ {\rm
Trace}\, (H_D^{-\frac{p}{2}}) \Bigr]^{1-\frac{p}{2}}\, .
\end{split}
\end{equation}
\end{proof}
In order to proceed further, we need the following intermediate result.
\begin{lem}
Let $u\in L^1_{\rm loc}([1,+\infty ))$ be a positive, locally
integrable function such that $u\in L^s ([1,+\infty ))$ for $s\in
(1,2]$ and
\begin{equation}
\int_1^\infty u(t)^s\, dt \le \frac{c}{s-1}\, \qquad s\in (1,2]
\end{equation}
for some constant $c>0$. Then there exists a constant $c^\prime
>0$ such that
\begin{equation}
\int_1^x u(t)\, dt\le c^\prime\ln x\, \qquad x\in (1,+\infty]\, .
\end{equation}
\end{lem}
\begin{proof}
By H\"older inequality and for $x\ge 1$, $s\in (1,2]$, we have
\[
\int_1^x u(t)\, dt \le \Bigl(\int_1^x u(t)^s\, dt\Bigr)^{1/s}
\cdot (x-1)^{1-1/s}
 \le \Bigl(\frac{c}{s-1}\Bigr)^{1/s} \cdot (x-1)^{1-1/s}\, .
\]
Setting $h(s):=\frac{1}{s}\ln\frac{c}{s-1} + \frac{s-1}{s}\ln
(x-1)$ we have $\int_1^x u(t)\, dt \le e^{h(s)}$. Evaluating
$h(s)$ at its critical point, where $\ln (x-1) = \frac{s}{s-1} +
\ln\frac{c}{s-1}$, we get $h(s)=1+\ln\frac{c}{s-1}$ and
\[
\int_1^x u(t)\, dt\le\frac{ec}{s-1}\, .
\]
As $s\in (1,2]$, we have $\ln (x-1) = \frac{s}{s-1} +
\ln\frac{c}{s-1} \ge \ln ec + \frac{1}{s-1}\ge\ln ce^2$, which
implies $x\ge  1 +ce^2$ and finally
\[
\int_1^x u(t)\, dt\le\frac{ec}{s-1} \le ec\ln \frac{x-1}{ec}\le
c^\prime \ln x
\]
for all $x\ge 1 +ce^2$ and $c^\prime := \alpha ec$ where $\alpha
>1$ is such that $ec\ge (\alpha -1 )^{(\alpha -1
)}/\alpha^\alpha$.
\end{proof}
We can now prove the finest summability properties for the quantum derivative
of functions with finite energy on fractals.
\begin{thm}{\bf (Commutator and Interpolation ideals)}
Let $(D,\mathbf{r})$ be a fixed regular harmonic structure on a
p.c.f. self-similar structure $\bigl(K , \S , \{F_i :i\in
\S\}\bigr)$ and let $\mu$ be a self-similar Borel measure on $K$ with
weights $(\mu_1 \, ,\dots\, ,\mu_N)$ such that $r_i\mu_i<1$ for
$\in\S$.\par\noindent Then $(F,\H)$ is a densely $(d_S ,
\infty)$-summable Fredholm module over $C(K)$:
\begin{equation}
[F,a]\in \mathcal{L}^{(d_S ,\infty )} (\H)\, \qquad a\in \F\,
\end{equation}
where $\mathcal{L}^{(d_S ,\infty )} (\H)$ is the interpolation ideal defined, for instance, in [C2 Chapter IV].
\end{thm}
\begin{proof}
By the upper bound (2.11) on the eigenvalue counting function,
there exists a constant $c_3 >0$ such that
\begin{equation}
\rho_* (x ,\mu )\le c_3 x^{d_S/2}, \qquad x\ge \lambda_1\, .
\end{equation}
As $k\le \rho_* (\lambda_k ,\mu )\le c_3 \lambda_k^{d_S/2}$, we
have $c_3^{-2/d_S} k^{2/d_S}\le \lambda_k$ and also
\begin{equation}
{\rm Trace} \Bigl(H_D^{-\frac{p}{2}}\Bigr) = \sum_{k=1}^\infty
\lambda_k^{-\frac{p}{2}}  \le
c_3^{\frac{p}{d_S}}\frac{d_S}{p-d_S}\,.
\end{equation}
As $p-d_S<1$, we have, for the constant $c_2 (p)$ in (3.14) the
bound
\begin{equation}
c_2 (p)\le 16c_1\Bigl(\frac{1}{\lambda_1} + 2 \Bigr) \Gamma
(p/2)^{-1}\frac{1}{p-d_S}\, .
\end{equation}
Combining (3.16), (3.29) and (3.30), we then have
\begin{equation}
{\rm Trace} \Bigl(|[F,a]|^p\Bigr) \le
\Bigl[16c_1\Bigl(\frac{1}{\lambda_1} + 2 \Bigr) \Gamma
(p/2)^{-1}\Bigr]^{\frac{p}{2}} c_3^{\frac{p}{d_S}(1-\frac{p}{2})}
d_S^{1-\frac{p}{2}}\frac{1}{p-d_S}
\end{equation}
so that for a suitable $c>0$ independent on $p\in (1,2]$
\begin{equation}
\sum_{k=1}^\infty \mu_k (T)^p \le c\frac{d_S}{p-d_S}\, ,\qquad d_S
< p\le 2\, ,
\end{equation}
where now $T:=[F,a]$. Setting $s:=p/d_S$ and $u(t) := \mu_{[t]}
(T)^{d_S}$ for $t\ge 1$, we have $s\in (1,2]$ and the thesis
follows applying the previous lemma:
\begin{equation}
\sup_{N\ge 2}\frac{1}{\ln N}\sum_{k=1}^N \mu_k (T)^{d_S} <+\infty
\end{equation}
so that $[F,a]\in \mathcal{L}^{(d_S ,\infty )} (\H)$ for all $a\in
\F$ as promised.
\end{proof}
Our final goal in this section is to provide a bound similar to
(3.16) in Theorem 3.6 but now involving Dixmier traces.
\begin{thm}{\bf (Dixmier trace summability)}
Let $(D,\mathbf{r})$ be a fixed regular harmonic structure on a
p.c.f. self-similar structure $\bigl(K , \S , \{F_i :i\in
\S\}\bigr)$, let $\mu$ be a self-similar Borel measure on $K$ with
weights $(\mu_1 \, ,\dots\, ,\mu_N)$ such that $r_i\mu_i<1$ for
$\in\S$ and $(F,\H)$ the associated densely $(d_S ,
\infty)$-summable Fredholm module over $C(K)$.\par\noindent Then,
for any Dixmier trace $\tau_\omega$, the following upper bound
holds true:
\begin{equation}
\tau_\omega \Bigl(|[F,a]|^{d_S}\Bigr) \le c_2
(d_S)^{\frac{d_S}{2}}
\cdot\Bigl(\E[a]\Bigr)^{\frac{d_S}{2}}\cdot\Bigl[\tau_\omega
\Bigl(H_D^{-\frac{d_S}{2}}\Bigl)\Bigr]^{1-\frac{d_S}{2}}\, ,\qquad\forall\,a\in
\F
\end{equation}
where $c_2 (d_S) := 32c_1 \Gamma (d_S /2)^{-1}$.
\end{thm}
\begin{proof}
By Theorem 3.8, $\tau_\omega \Bigl(|[F,a]|^{d_S}\Bigr)$ is finite
for all Dixmier functionals $\omega$ on $L^\infty
(\mathbb{R}^*_+)$ and, by [CPS Lemma 5.1], we have the identity
\begin{equation}
d_S\tau_\omega (|[F,a]|^{d_S}) = \widetilde\omega -
\lim_{r\to\infty}\frac{1}{r} \tau (|[F,a]|^{d_S+\frac{1}{r}})
\end{equation}
where $\widetilde\omega :=\omega\circ L$ is the Dixmier functional
on $L^\infty (\mathbb{R})$ corresponding to $\omega$ through the
map $L:L^\infty (\mathbb{R})\to L^\infty (\mathbb{R}^*_+)$ given
by $Lf := f\circ\log$.\par\noindent By Lemma 3.7 applied to the
bound (3.29), we have that $\tau_\omega
\Bigl(H_D^{-\frac{d_S}{2}}\Bigr)$ is finite for all Dixmier
functionals $\omega$ on $L^\infty (\mathbb{R}^*_+)$ so that, again
by [CPS Lemma 5.1], we have the identity
\begin{equation}
d_S\tau_\omega \Bigl(H_D^{-\frac{d_S}{2}}\Bigr) = \widetilde\omega
- \lim_{r\to\infty}\frac{1}{r} \tau\omega
\Bigl(H_D^{-\frac{d_S+\frac{1}{r}}{2}}\Bigr)\, .
\end{equation}
The desired bound (3.34) then follows by (3.16) in Theorem 3.6.
\end{proof}
The previous result naturally suggests the consideration of a new energy functional
which should be a conformal invariant in the sense of Alain Connes [C2].
\begin{defn}
The functional $\Phi^{d_S}_\omega : \F_0 \longrightarrow [0,+\infty
)$
\begin{equation}
\Phi^{d_S}_\omega (a) := \tau_\omega \Bigl(|[F,a]|^{d_S}\Bigr)
\end{equation}
will be referred to as the $d_S$-energy functional of the harmonic
structure $(D,\mathbf{r})$.
\end{defn}
\begin{cor}
For all $a\in\F$ we have
\begin{equation}
\sum_{i=1}^N \Phi^{d_S}_\omega (a\circ F_i ) \le c_2
(d_S)^{\frac{d_S}{2}}
\cdot\Bigl(\E[a]\Bigr)^{\frac{d_S}{2}}\cdot\Bigl[\tau_\omega
\Bigl(H_D^{-\frac{d_S}{2}}\Bigl)\Bigr]^{1-\frac{d_S}{2}} \, .
\end{equation}
\end{cor}
\begin{proof}
Setting $c:=c_2 (d_S)^{\frac{d_S}{2}} \cdot\Bigl[\tau_\omega
\Bigl(H_D^{-\frac{d_S}{2}}\Bigl)\Bigr]^{1-\frac{d_S}{2}}$ and, for
$a\in\F$, applying (3.34) to $a\circ F_i\in\F$, we have
\begin{equation}
\Phi^{d_S}_\omega (a\circ F_i ) \le c \cdot\Bigl(\E[a\circ F_i
]\Bigr)^{\frac{d_S}{2}}\, .
\end{equation}
By H\"older inequality we then have
\[
\begin{split}
\sum_{i=1}^N \Phi^{d_S}_\omega (a\circ F_i ) &\le c\cdot
\sum_{i=1}^N \Bigl(\E[a\circ F_i ]\Bigr)^{\frac{d_S}{2}}\\
& = c\cdot
\sum_{i=1}^N r_i^{\frac{d_S}{2}}\Bigl(r_i^{-\frac{d_S}{2}}\E[a\circ F_i ]\Bigr)^{\frac{d_S}{2}}\\
& \le c\cdot \Bigl(\sum_{i=1}^N
r_i^{\frac{d_S}{2}\frac{2}{2-d_S}}\Bigl)^{\frac{2-d_S}{2}}\cdot
\Bigl(\sum_{i=1}^Nr_i^{-1}\E[a\circ F_i ]\Bigr)^{\frac{d_S}{2}}\\
& = c\cdot \Bigl(\sum_{i=1}^N
r_i^{d_H}\Bigl)^{\frac{2-d_S}{2}}\cdot
\Bigl(\E[a ]\Bigr)^{\frac{d_S}{2}}\\
& = c\cdot
\Bigl(\E[a ]\Bigr)^{\frac{d_S}{2}}\\
\end{split}
\]
\end{proof}
The previous result suggests that the $d_S$-energy functional may
be conformal, as we now prove that it is indeed, by means of the uniqueness result of [CS].
\begin{thm}{\bf (Conformal invariance)}
The $d_S$-energy functional is a self similar conformal invariant
\begin{equation}
\Phi^{d_S}_\omega (a)= \sum_{i=1}^N \Phi^{d_S}_\omega (a\circ F_i
)\qquad a\in \F\, .
\end{equation}
\end{thm}
\begin{proof}
Let us consider the Hilbert space $\H^N = \bigoplus_{i=1}^N \H$
endowed with the action of $C(K)$ given by
\[
a\Bigl(\bigoplus_{i=1}^N \xi_i\Bigr) := \bigoplus_{i=1}^N (a\circ
F_i)\xi_i\, ,\qquad a\in C(K)\, ,\qquad \xi_i\in \H\, ,\qquad
i=1,\cdots ,N\, ,
\]
and the involution given by $J^N :\H^N\to \H^N$
\[
J^N \Bigl(\bigoplus_{i=1}^N \xi_i\Bigr) := \bigoplus_{i=1}^N
J\xi_i\, ,\qquad \xi_i\in \H\, .
\]
It is easily verified that $(C(K), \H^N , J^N )$ is a symmetric
Hilbert module over $C(K)$ and the map $\partial^N :\F\to\H^N$
given by
\[
\partial^N (a):=\bigoplus_{i=1}^N r_i^{-1/2}\partial (a\circ F_i)
\]
is a symmetric derivation such that
\[
\|\partial^N (a)\|^2_{\H^N} = \bigoplus_{i=1}^N r_i^{-1}
\|\partial (a\circ F_i)\|^2_{\H} = \bigoplus_{i=1}^N r_i^{-1} \E
[a\circ F_i ] = \E [a]\, .
\]
In other words, $(\partial^N\, ,\H^N , J^N )$ is a new symmetric
derivation representing the Dirichlet form $(\E ,\F )$, isomorphic
to the older one $(\partial\, ,\H , J )$ by [CS Theorem 8.3].
Since the corresponding Fredholm modules are unitarily isomorphic,
the $d_S$-energy functional $\Phi^{d_S}_\omega$ is unchanged if
computed using the new structure.

\end{proof}

\end{document}